\documentclass[a4paper,12pt]{amsart}
\usepackage{mathrsfs}
\usepackage{amsfonts}
\usepackage{amssymb}
\usepackage{ifthen}
\usepackage{graphicx}
\usepackage[usenames]{color}

\nonstopmode \numberwithin{equation}{section}
\newtheorem{thm}{Theorem}[section]
\newtheorem{lem}{Lemma}[section]
\newtheorem{cor}{Corollary}[section]

\newtheorem{conj}[equation]{Conjecture}

\theoremstyle{definition}
\newtheorem{defn}{Definition}[section]

\newtheorem{prob}[equation]{Problem}
\newtheorem{ques}[equation]{Question}
\newtheorem{rem}{Remark}[section]
\newtheorem{exam}[equation]{Example}

\newcounter {own}
\def\theown {\thesection       .\arabic{own}}

\newenvironment{pf}[1][]{%
 \vskip 3mm
 \noindent
 \ifthenelse{\equal{#1}{}}%
  {{\slshape Proof. }}%
  {{\slshape #1.} }%
 }%
{\qed\bigskip}

\newcounter{alphabet}
\newcounter{tmp}
\newenvironment{Thm}[1][]{\refstepcounter{alphabet}%
\bigskip%
\noindent%
{\bf Theorem \Alph{alphabet}}%
\ifthenelse{\equal{#1}{}}{}{ (#1)}%
{\bf .} \itshape}{\vskip 8pt}

\makeatletter
\newcommand{\Ref}[1]{\@ifundefined{r@#1}{}{\setcounter{tmp}{\ref{#1}}\Alph{tmp}}}
\makeatother

\newenvironment{Lem}[1][]{\refstepcounter{alphabet}%
\bigskip%
\noindent%
{\bf Lemma \Alph{alphabet}}%
{\bf .} \itshape}{\vskip 8pt}

\newcommand{\ID}{{\mathbb D}}

\def\be{\begin{equation}}
\def\ee{\end{equation}}

\newcommand{\bee}{\begin{enumerate}}
\newcommand{\eee}{\end{enumerate}}

\newcommand{\blem}{\begin{lem}}
\newcommand{\elem}{\end{lem}}
\newcommand{\bthm}{\begin{thm}}
\newcommand{\ethm}{\end{thm}}
\newcommand{\bcor}{\begin{cor}}
\newcommand{\ecor}{\end{cor}}
\newcommand{\beg}{\begin{exam}}
\newcommand{\eeg}{\end{exam}}
\newcommand{\begs}{\begin{examples}}
\newcommand{\eegs}{\end{examples}}
\newcommand{\bdefe}{\begin{defn}}
\newcommand{\edefe}{\end{defn}}
\newcommand{\bprob}{\begin{prob}}
\newcommand{\eprob}{\end{prob}}
\newcommand{\bques}{\begin{ques}}
\newcommand{\eques}{\end{ques}}
\newcommand{\bei}{\begin{itemize}}
\newcommand{\eei}{\end{itemize}}
\newcommand{\bcon}{\begin{conj}}
\newcommand{\econ}{\end{conj}}
\newcommand{\bcons}{\begin{conjs}}
\newcommand{\econs}{\end{conjs}}
\newcommand{\bprop}{\begin{propo}}
\newcommand{\eprop}{\end{propo}}
\newcommand{\br}{\begin{rem}}
\newcommand{\er}{\end{rem}}
\newcommand{\brs}{\begin{rems}}
\newcommand{\ers}{\end{rems}}
\newcommand{\bo}{\begin{obser}}
\newcommand{\eo}{\end{obser}}
\newcommand{\bos}{\begin{obsers}}
\newcommand{\eos}{\end{obsers}}
\newcommand{\bpf}{\begin{pf}}
\newcommand{\epf}{\end{pf}}
\newcommand{\ba}{\begin{array}}
\newcommand{\ea}{\end{array}}
\newcommand{\beq}{\begin{eqnarray}}
\newcommand{\beqq}{\begin{eqnarray*}}
\newcommand{\eeq}{\end{eqnarray}}
\newcommand{\eeqq}{\end{eqnarray*}}

\newcommand{\ds}{\displaystyle}

\begin{document}
\bibliographystyle{amsplain}
\title {Uniformly locally univalent harmonic mappings}

\author{Saminathan Ponnusamy 
}
\address{S. Ponnusamy,
Indian Statistical Institute (ISI), Chennai Centre, SETS (Society
for Electronic Transactions and Security), MGR Knowledge City, CIT
Campus, Taramani, Chennai 600 113, India.}
\email{samy@isichennai.res.in, samy@iitm.ac.in}
\author{Jinjing  Qiao $^\dagger $}
\address{J. Qiao, Department of Mathematics,
Hebei University, Baoding, Hebei 071002, People's Republic of China.}
\email{mathqiao@126.com}
\author{Xiantao Wang}
\address{X. Wang, Department of Mathematics, Shantou University, Shantou, Guangdong  515063,
People's Republic of China.}
\email{xtwang@stu.edu.cn}

\keywords{Harmonic mapping, Pre-Schwarzian derivatives, uniformly locally univalence, growth estimate, coefficient estimate, harmonic Bloch space,
Hardy space.\\
${} ^\dagger$ Corresponding author.
}

\begin{abstract}
The primary aim of this paper is to characterize the uniformly locally univalent
harmonic mappings in the unit disk. Then, we obtain
sharp distortion, growth and covering theorems for one parameter family
${\mathcal B}_{H}(\lambda)$ of uniformly locally univalent
harmonic mappings. Finally, we show that the subclass of
$k$-quasiconformal harmonic mappings in ${\mathcal B}_{H}(\lambda)$
and the class ${\mathcal B}_{H}(\lambda)$ are
contained in the Hardy space of a specific exponent depending on the
$\lambda$, respectively,  and we also discuss the growth of
coefficients for harmonic mappings in ${\mathcal B}_{H}(\lambda)$.
\end{abstract}

\thanks{}

\maketitle
\pagestyle{myheadings}
\markboth{S. Ponnusamy, J. Qiao and X. Wang}{Uniformly locally univalent harmonic mappings}

\section{Introduction}\label{csw-sec1}

The class of complex-valued harmonic mappings $f$ defined on a
simply connected domain $D$ of the complex plane $\mathbb{C}$ has attracted
the attention of function theorists because it
generalizes the class of analytic functions with a lot of rich applications
in many different fields. Every such $f$ has the canonical decomposition $f=h+\overline{g}$,
where both $h$ and $g$ are analytic in $D$ and $g(z_0)=0$ for some prescribed point $z_0\in D$
(cf. \cite{CS,Du}). For a complex-valued and continuously differentiable mapping $f$, we let
$$\lambda_{f}=|f_{z}|-|f_{\overline{z}}|  ~\mbox{ and }~
\Lambda_{f}=|f_{z}|+|f_{\overline{z}}|
$$
so that the Jacobian $J_{f}$ of $f$ takes the form
$$J_{f}=\lambda_{f}\Lambda_{f}=|f_{z}|^{2}-|f_{\overline{z}}|^{2}.
$$
Moreover, a necessary and sufficient condition for harmonic mappings
$f=h+\overline{g}$ to be locally univalent and sense preserving in
$D$ is that $J_f=|h'|^2-|g'|^2>0$, or equivalently, its dilatation
the dilatation $\omega_f(z)=g'(z)/h'(z)$ satisfies the inequality $|\omega_f(z)|<1$ for $z\in D$.
(see \cite{Le} and \cite{CS,Du,PonRasi2013}).
Let $k\in [0,1)$ be a constant. Then, we say that a sense preserving harmonic
mapping $f=h+\overline{g}$ in $D$ is a $k$-quasiconformal mapping
if $|\omega_f(z)|\leq k$ holds in $D$.


Throughout this paper, we consider harmonic mappings defined on the unit disk $\mathbb{D}=\{z\in \mathbb{C}:\, |z|<1\}$.
Denote by  ${\mathcal H}$ the class of harmonic mappings $f=h+\overline{g}$ in $\mathbb{D}$
such that $h(0)=g(0)=h'(0)-1=0$ and consider the family
$${\mathcal S}_{H}=\{f\in {\mathcal H}:\, f \mbox{ is sense-preserving and univalent in $\mathbb{D}$}
\}.
$$
Often it is convenient to work with
$${\mathcal S}_{H}^0=\{f\in {\mathcal S}_{H}:\, f_{\overline{z}}(0)=0\}.
$$
Although  both the families ${\mathcal S}_{H}$ and ${\mathcal S}_{H}^0$ are known to be normal, only
${\mathcal S}_{H}^0$ is compact (see \cite{CS}). We also denote the class of analytic functions
$f$ in $\mathbb{D}$  with $f(0)=f'(0)-1=0$ by ${\mathcal A}$ so that $\mathcal H$ reduces to
${\mathcal A}$ when the co-analytic part $g$ of $f=h+\overline{g}\in\mathcal H$ vanishes identically in $\ID$.
Then the set ${\mathcal S}:={\mathcal A}\cap {\mathcal S}_{H}^0$ of all
normalized univalent analytic functions in $\ID$
is the central object in the study of geometric function theory so that
$ {\mathcal S}\subset {\mathcal S}_{H}^0\subset {\mathcal S}_{H}.$

We denote $d_h(z,w)$ as the hyperbolic distance of $z, w \in \mathbb{D}$, that is,
$$d_h(z,w)=\frac{1}{2}\log\left (\frac{1+\left |\frac{z-w}{1-\overline{z}w}\right |}{1-\left |\frac{z-w}{1-\overline{z}w}\right |}\right ).
$$
We say that a harmonic mapping
$f=h+\overline{g}$ in $\mathbb{D}$ is uniformly locally univalent if
$f$ is univalent in each hyperbolic disk
$$D_h(a,\rho)=\left \{z\in \mathbb{D}:\, d_h(z,a)< \rho\right \}
$$
with center $a\in \mathbb{D}$ and hyperbolic radius $\rho$ (independent of the center), $0<\rho\leq \infty$.
The subscript $h$ in $d_h$ and $D_h$ should not be confused with the analytic part $h$ of the harmonic mapping $f$.

If $f$ is analytic in the above definition, then it reduces to the uniformly
locally univalent (analytic) functions. We know that a holomorphic
universal covering map of a plane domain $D$ is uniformly locally
univalent if and only if the boundary of $D$ is uniformly perfect
(cf. \cite{Ch, Su}). Also, in \cite{Ya}, Yamashita  showed that an
analytic function $f$ in $\mathbb{D}$ is uniformly locally univalent in $\ID$
if and only if the pre-Schwarzian derivative  $T_f=f''/f'$ of $f$ is
hyperbolically bounded, i.e., the norm
$$\|T_f\|=\sup_{z\in \mathbb{D}}(1-|z|^2)|T_f(z)|
$$
is finite and this means that $\log f'$ belongs to the Bloch space
$\mathscr{B}$ (cf. \cite{ACP, Co1}).

In Section \ref{sec-2} (see Theorem \ref{thm2.1}), we characterize the
uniformly locally univalent harmonic mappings $f=h+\overline{g}$ in
terms of the pre-Schwarzian derivative of $h+e^{i\theta}g$ for each $\theta\in [0,2\pi]$.
This result and the corresponding results in \cite{KS} helps to
obtain  sharp distortion, growth and covering theorems (see Section
\ref{sec-3}) for the class ${\mathcal B}_{H}(\lambda)$,
where $\lambda$ is a positive real number, and
$${\mathcal B}_{H}(\lambda)=\{f=h+\overline{g}\in {\mathcal H}:\, \|T_f\|\leq 2\lambda\}
$$ with
\be\label{eq1.1} \|T_f\|:=\sup_{z\in \mathbb{D},\ \theta\in [0,
2\pi]} (1-|z|^2)\left
|\frac{h''(z)+e^{i\theta}g''(z)}{h'(z)+e^{i\theta}g'(z)}\right|.
\ee
Henceforth, $\|T_f\|$ is defined by \eqref{eq1.1} in the case of harmonic mappings $f=h+\overline{g}$ in$\ID$.


It is known that for $\lambda>1$, the class ${\mathcal B}(\lambda)$  and
the subclass ${\mathcal B}(\lambda)\cap{\mathcal S}$ are contained in the Hardy space
$H^p$ with $0<p<1/(\lambda^2-1)$ and $0<p<1/(\lambda-1)$, respectively (cf. \cite{Ki, KS1}).

In Section \ref{sec-4'},  we consider relationships between the
 space ${\mathcal B}_{H}(\lambda)$ and the harmonic Hardy space.   We
also prove that a $k$-quasiconformal harmonic mapping $f\in
{\mathcal B}_{H}(\lambda)$ ($\lambda>1$) is contained in the
harmonic Hardy space $h^p$ with $0<p<1/(\lambda-1)$, and also obtain
that ${\mathcal B}_{H}(\lambda)\subset h^p$ with
$0<p<1/(\lambda^2-1)$. Finally, in the last section, as applications
of distortion estimate obtained in Section \ref{sec-3}, we discuss
the growth of coefficients for harmonic mappings in ${\mathcal
B}_{H}(\lambda)$ ($\lambda>1)$.

In \cite{KS}, the authors discussed the set ${\mathcal B}(\lambda):={\mathcal A}\cap {\mathcal
B}_{H}(\lambda)$ and obtained distortion estimates for analytic
functions in  ${\mathcal B}(\lambda)$ in terms of $\lambda$, and
characterization for functions in ${\mathcal B}(\lambda)$ (cf.
\cite[Proposition 1.1]{KS}). As a consequence of  Theorem
\ref{thm2.1} in Section \ref{sec-2} and \cite[Proposition 1.1]{KS},
we can easily obtain the following corollary which characterizes
harmonic mappings in ${\mathcal B}_{H}(\lambda)$. We omit its proof
and this particular case is indeed a generalization of earlier known
result (see \cite[Proposition 1.1]{KS}) to the case of harmonic
mappings.

\bcor \label{thm2.2}
A locally univalent harmonic mapping $f=h+\overline{g}\in {\mathcal H}$
belongs to ${\mathcal B}_{H}(\lambda)$ if and only if,  for each
pair of points $z_1$, $z_2$ in $\mathbb{D}$ and $\theta\in
[0,2\pi]$,
$$|u_\theta(z_1)-u_\theta(z_2)|\leq 2\lambda d_h(z_1, z_2),
$$
where $u_\theta(z)=\log \big(h'(z)+e^{i\theta}g'(z)\big)$.
\ecor

\section{Characterizations of Uniformly locally univalent harmonic mappings}\label{sec-2}

We now state our first result which is indeed a generalization of \cite[Theorem 1]{Ya} to the case of harmonic
mappings.

\bthm\label{thm2.1} A harmonic mapping $f=h+\overline{g}$ is
uniformly locally univalent in $\mathbb{D}$ if and only if $\|T_f\|<\infty$.
\ethm


For the proof of the sufficiency of Theorem \ref{thm2.1}, we need the following classical result due
to Noshiro \cite{No}.

\begin{Lem}\label{Noshi-lemA}
Let $f(z)=z+\sum_{k=2}^{\infty}a_{k}z^{k}$ be analytic for $|z|<R$ and $|f'(z)|<M$ for $|z|<R$.
Then the disk $|z|<R/M$ is mapped on a starlike domain with respect to the origin by $f$
and also by all its polynomial sections $f_n(z)=z+\sum_{k=2}^{n}a_{k}z^{k}$ $~(n=2,3, \ldots)$.
\end{Lem}

\subsection*{Proof of Theorem \ref{thm2.1}}
Let $f=h+\overline{g}$ be harmonic in $\mathbb{D}$ and assume that $\|T_f\|<\infty$.
Define $F(\xi)=(f\circ T)(\xi)$ for $\xi\in \mathbb{D}$, where
$$w=T(\xi)=\frac{R\xi+a}{1+\overline{a}R \xi}
$$
with $R=\tanh \rho$, the constants $a\in \mathbb{D}$ and $\rho$ ($0<\rho\leq \infty$). Then $F=H+\overline{G}$ is harmonic in $\mathbb{D}$.
 Elementary computations yield \beqq
\frac{H''(\xi)+e^{i\theta}G''(\xi)}{H'(\xi)+e^{i\theta}G'(\xi)}
&=&\frac{h''(w)+e^{i\theta}g''(w)}{h'(w)+e^{i\theta}g'(w)}T'(\xi)
+\frac{T''(\xi)}{T'(\xi)},
\eeqq
where
$$T'(\xi)=\frac{R(1-|a|^2)}{(1+\overline{a}R \xi)^2}~\mbox{ and }~
T''(\xi)=-\frac{2\overline{a}(1-|a|^2)R ^2}{(1+\overline{a}R \xi)^3}
$$
so that
$$\frac{T''(\xi)}{T'(\xi)}=-\frac{2\overline{a}R}{1+\overline{a}R \xi}.
$$
Since
$$1-|w|^2=\frac{(1-|a|^2)(1-|\xi|^2 R ^2 )}{|1+\overline{a}R \xi|^2},
$$
we easily have $|T'(\xi)|(1-|\xi|^2)\leq 1-|w|^2$ and therefore, it follows that
\beqq
(1-|\xi|^2)\left |\frac{H''(\xi)+e^{i\theta}G''(\xi)}{H'(\xi)+e^{i\theta}G'(\xi)}\right|
&\leq &
(1-|w|^2)\left |\frac{h''(w)+e^{i\theta}g''(w)}{h'(w)+e^{i\theta}g'(w)}\right|\\
&& \hspace{.2cm} +(1-|\xi|^2)\left
|\frac{2\overline{a}R}{1+\overline{a}R \xi}\right |. \eeqq This
inequality implies that \be\label{eq2.2} \sup_{\xi\in
\mathbb{D}}(1-|\xi|^2)\left
|\frac{H''(\xi)+e^{i\theta}G''(\xi)}{H'(\xi)+e^{i\theta}G'(\xi)}\right
| \leq  k_0<\infty , \ee where
$$ k_0=\sup_{w \in \mathbb{D}}(1-|w|^2)\left|\frac{h''(w)+e^{i\theta}g''(w)}{h'(w)+e^{i\theta}g'(w)}\right |
+\frac{2R}{1-R}.
$$

Let $\varphi$ be an analytic branch of $\log\big(H'(\xi)+e^{i\theta}G'(\xi)\big)$ in $\mathbb{D}$. Then
$$\varphi'(\xi)=\frac{H''(\xi)+e^{i\theta}G''(\xi)}{H'(\xi)+e^{i\theta}G'(\xi)}.
$$
This choice is clearly possible, because $H'(\xi)+e^{i\theta}G'(\xi)\neq 0$
for $\xi\in \mathbb{D}$, by \eqref{eq2.2}. It then follows from \eqref{eq2.2} that
\be\label{eq2.3}
\left|\log\Big|\frac{H'(\xi)+e^{i\theta}G'(\xi)}{H'(0)+e^{i\theta}G'(0)}\Big|\right|\leq
|\varphi(\xi)-\varphi(0)|\leq
\frac{k_0}{2}\log\left(\frac{1+|\xi|}{1-|\xi|}\right).
\ee
Now, we introduce $H_{\theta}(\xi)$ by
$$H_{\theta}(\xi)=\frac{H(\xi)+e^{i\theta}G(\xi)}{H'(0)+e^{i\theta}G'(0)}.
$$
We see that $H_{\theta}$ is analytic in $\ID$ and is normalized so
that $H_{\theta}'(0)-1=0$. We infer from \eqref{eq2.3} that
$$\log \big|H_{\theta}'(\xi)\big|\leq
\frac{k_0}{2}\log3 ~\mbox{ for }~ |\xi|< \frac{1}{2},
$$
whence
$$\big|H_{\theta}'(\xi)\big|<3^{k_0/2} ~\mbox{ for }~ |\xi|< \frac{1}{2}.
$$
Therefore, by  Lemma \Ref{Noshi-lemA} of Noshiro,  $H_{\theta}(\xi)-H_{\theta}(0)$ is
univalent in the disk $|\xi|<\frac{3^{-k_0/2}}{2}$ for each $\theta$. The radius of convexity
for univalent functions is known to be $2-\sqrt{3}$
(cf. \cite[Theorem 2.13]{Du}). Thus, $H_{\theta}(\xi)-H_{\theta}(0)$
and $H(\xi)+e^{i\theta}G(\xi)$ are convex in
$|\xi|<(2-\sqrt{3})\frac{3^{-k_0/2}}{2}=\rho_0$. This implies that
$F$ is harmonic (convex) univalent in
$|\xi|<(2-\sqrt{3})\frac{3^{-k_0/2}}{2}$ (cf. \cite{CS}).

Consequently, $f$ is univalent in the hyperbolic disk $D_h(a, \rho_0)$ with
$\tanh\rho_0=(2-\sqrt{3})\frac{3^{-k_0/2}}{2}\tanh\rho$ if
$(2-\sqrt{3})\frac{3^{-k_0/2}}{2}\leq 1$, and $\rho_0=\rho$ if
$(2-\sqrt{3})\frac{3^{-k_0/2}}{2}> 1$. Hence, $f$ is uniformly
locally univalent.

To prove the necessity, we assume that $f$ is uniformly locally
univalent in $\mathbb{D}$, that is,  $f$ is univalent in each hyperbolic disk $D_h(a,\rho)$,
where $a\in \mathbb{D}$ and $0<\rho\leq \infty$.

Again, as above, we consider $w=T(\xi)$ and let
$$F(\xi)=(f\circ T)(\xi)=f(w) ~\mbox{ for $\xi\in \mathbb{D}$}.
$$
Then $F=H+\overline{G}$ is univalent in $\mathbb{D}$. By Lewy's Theorem (cf. \cite{Du}),
the Jacobian $J_F$ is different from $0$ for all $z\in\mathbb{D}$,
and then, without loss of generality, we assume that $F$ is sense-preserving. Let
$$F_0(\xi)=\frac{H(\xi)-H(0)}{H'(0)}+\overline{\frac{G(\xi)-G(0)}{H'(0)}}=H_0(\xi)+\overline{G_0(\xi)}.
$$
Obviously, $F_0\in {\mathcal S}_{H}$. For $\xi\in \mathbb{D}$, set
$$F_1(z)=\frac{F_0\Big(\frac{z+\xi}{1+\overline{\xi}z}\Big)-F_0(\xi)}{(1-|\xi|^2)H_0'(\xi)}
= H_1(z)+\overline{G_1(z)},
$$
which again belongs to ${\mathcal S}_{H}$. The analytic function $H_1(z)$ has the form
$$H_1(z)=z+A_2(\xi)z^2+A_3(\xi)z^3+\cdots
$$
and a direct computation shows that
$$A_2(\xi)=\frac{1}{2}\left \{ (1-|\xi|^2)\frac{H_0''(\xi)}{H_0'(\xi)}-2\overline{\xi} \right \}
=\frac{1}{2}\left \{ (1-|\xi|^2)\frac{H''(\xi)}{H'(\xi)}-2\overline{\xi} \right \}.
$$
Let $\alpha=\sup\{|a_2|:\,
f(z)=\sum_{k=1}^{\infty}a_kz^k+\sum_{k=1}^{\infty}\overline{b}_k\overline{z}^k\in{\mathcal
S}_H \}$. For $f\in \mathcal {S}_H$, we have
$\frac{f(z)-\overline{b_1f(z)}}{1-|b_1|^2}\in \mathcal {S}_H^0$. It
is known that for
$f^*(z)=\sum_{k=1}^{\infty}a^*_kz^k+\sum_{k=1}^{\infty}\overline{b}^*_k\overline{z}^k\in
\mathcal {S}_H^0$, the coefficient $|a^*_2|<49$ and
$|b^*_2|<\frac{1}{2}$(cf. \cite{Du}). Using this estimate, by  computations, it is possible to get
$|a_2|<98$. It has been recently shown by Abu-Muhanna et al. \cite{AAP2015} that  $|a^*_2|\leq 16.5$
which indeed the best known upper bound for $|a^*_2|$.  Since $F_1\in {\mathcal S}_{H}$, we must have
$|A_2(\xi)|\leq \alpha$ and therefore,
$$(1-|\xi|^2)\left |\frac{H''(\xi)}{H'(\xi)}\right |< 2(\alpha+1), ~\mbox{ $\xi\in \mathbb{D}$}.
$$

For each $c\in \mathbb{D}$, the composition of sense preserving affine mapping $\phi(w) =w+c\overline{w}$
with $F$, namely, the function $F+c\overline{F}$, is univalent and sense-preserving in $\ID$. Then by what we have just proved, we obtain
$$(1-|\xi|^2)\left |\frac{H''(\xi)+cG''(\xi)}{H'(\xi)+cG'(\xi)}\right |< 2(\alpha+1), ~\mbox{ $\xi\in \mathbb{D}$},
$$
which in particular implies that, for each $\theta\in [0, 2\pi]$,
$$ (1-|\xi|^2)\left |\frac{H''(\xi)+e^{i\theta}G''(\xi)}{H'(\xi)+e^{i\theta}G'(\xi)}\right |
< 2(\alpha+1), ~\mbox{ $\xi\in \mathbb{D}$}.
$$
Thus, for $f=h+\overline{g}$, we have
$$ A(\theta):=\sup_{z\in \mathbb{D}}
(1-|z|^2)\left |\frac{h''(z)+e^{i\theta}g''(z)}{h'(z)+e^{i\theta}g'(z)}\right |
<\infty.
$$
Since $A(\theta)$ is a continuous function of $\theta$ in $[0, 2\pi]$, it follows from $ A(\theta)<\infty$ that
$$ \sup_{z\in \mathbb{D},\ \theta\in [0, 2\pi]} (1-|z|^2)\left
|\frac{h''(z)+e^{i\theta}g''(z)}{h'(z)+e^{i\theta}g'(z)}\right
|<\infty.
$$
The proof of the theorem is complete.
\hfill $\Box$

\section{Growth estimate for the class ${\mathcal B}_{H}(\lambda)$}\label{sec-3}

For a nonnegative real number $\lambda$, we consider
$$H_\lambda(z)=\int_{0}^z\Big(\frac{1+t}{1-t}\Big)^{\lambda}\,dt.
$$
It is easy to verify that $\|T_{H_\lambda}\|=2\lambda$, and thus
$H_\lambda\in {\mathcal B}_H(\lambda)$. If $\lambda \geq 0$,  then
it is known that  $H_\lambda$ is univalent in $\mathbb{D}$ if and
only if $0\leq\lambda\leq1$ (see \cite[Lemma 2.1]{KS}). We will see
later that $H_\lambda$ is extremal in the class ${\mathcal
B}_{H}(\lambda)$. It follows from
Theorem \ref{thm2.1} that if $f=h+\overline{g}\in {\mathcal
B}_{H}(\lambda)$, then $\frac{h+e^{i\theta}g}{1+e^{i\theta}b_1}\in
{\mathcal B}(\lambda)$. This fact and \cite[Theroem
2.3]{KS} give the following result.

\bthm[Distortion theorem]\label{thm3.1}
Let $\lambda$ be a nonnegative real number and
$f(z)=h(z)+\overline{g(z)}=\sum_{n=1}^{\infty}a_nz^n+\sum_{n=1}^{\infty}\overline{b}_n\overline{z}^n\in
{\mathcal B}_{H}(\lambda)$. Then for $z\in \ID$, we have
$$ \big | \lambda_{f} (z)\big |=\big|\,|h'(z)|-|g'(z)|\, \big|\geq |1-|b_1|\,| \left (\frac{1-|z|}{1+|z|}\right )^\lambda=|1-|b_1|\,| H_\lambda'(-|z|),
$$
$$\big | \Lambda_{f}(z) \big |= |h'(z)|+|g'(z)|\leq (1+|b_1|) \left (\frac{1+|z|}{1-|z|}\right )^\lambda=(1+|b_1|) H_\lambda'(|z|)
$$
and $|f(z)|\leq (1+|b_1|)H_\lambda(|z|)$. Furthermore, if $f\in
\mathcal {S}_H^0\cap{\mathcal B}_{H}(\lambda)$, then
$$-H_\lambda(-|z|)\leq |f(z)|\leq H_\lambda(|z|).
$$
Equality occurs in each case when $f(z)=\overline{\mu}H_\lambda(\mu z)$ for a unimodular constant $\mu$.
\ethm

\bcor\label{cor3.1} For $\lambda>1$, each $f(z)=h(z)+\overline{g(z)}=\sum_{n=1}^{\infty}a_nz^n+\sum_{n=1}^{\infty}\overline{b}_n\overline{z}^n\in
{\mathcal B}_{H}(\lambda)$ satisfies the growth condition
$$f(z)=O\big((1-|z|)^{1-\lambda}\big)
$$
as $|z|\rightarrow 1$. On the other hand, for $\lambda<1$, each
mapping $f\in {\mathcal B}_{H}(\lambda)$ is bounded with the bound
$(1+|b_1|)H_\lambda(1)$. Moreover,  if  $\lambda>0$ and $f\in
\mathcal {S}_H^0\cap{\mathcal B}_{H}(\lambda)$ in $\ID$, then the
image $f(\mathbb{D})$ contains the disk $\{w:\,|w|<-H_\lambda(-1)\}$.
\ecor

By \cite{Be, BC}, for $\lambda\leq 1/2$, ${\mathcal
B}(\lambda)\subset {\mathcal S}$ and so, by Theorem \ref{thm2.1},
for $\lambda\leq 1/2$, $f\in {\mathcal B}_{H}(\lambda)$ must be
univalent in $\ID$. We also note that, for $0\leq \lambda\leq 1$, we
have
$$-H_\lambda(-1)\geq -H_1(-1)=2\log2-1=0.38629 \ldots,
$$
and therefore the result is an improvement of the covering theorem
for harmonic mappings in ${\mathcal S}_{H}^0$.

In Corollary \ref{cor3.1}, the case $\lambda=1$ is critical. By Theorem \ref{thm3.1}, we have that, for
$f\in {\mathcal B}_{H}(1)$,
$$|f(z)|\leq (1+|b_1|)H_1(|z|)=(1+|b_1|)\big(-2\log(1-|z|)-|z|\big),
$$
which shows that functions in ${\mathcal B}_{H}(1)$ need not be bounded. The next theorem, which follows from
Theorem \ref{thm2.1} and \cite[Proposition 2.5]{KS}, gives a boundedness criterion for mappings in ${\mathcal B}_{H}(1)$.

\bthm\label{thm3.2}
If a harmonic mapping $f=h+\overline{g}$ in $\mathbb{D}$ satisfies the condition
$$ {\overline{\lim_{|z|\rightarrow 1^{-}}}}
\left \{(1-|z|^2)\left |\frac{h''(z)+e^{i\theta}g''(z)}{h'(z)+e^{i\theta}g'(z)}\right |-2\right \}
\log\frac{1}{1-|z|^2}<-2
$$
for each $\theta\in [0,2\pi]$, then $f$ is bounded. Here the constant $-2$ on the right hand side is
sharp.
\ethm

We conclude this section with the H\"{o}lder continuity of mappings in ${\mathcal B}_{H}(\lambda)$.

\bthm\label{thm3.3}
Let $0\leq \lambda<1$. Then each mapping $f\in {\mathcal B}_{H}(\lambda)$ is
H\"{o}lder continuous of exponent $1-\lambda$ in $\mathbb{D}$.
\ethm

The proof follows from Theorem \ref{thm2.1} and \cite[Theorem 2.6]{KS} and so, we omit its detail.

\section{The space ${\mathcal B}_{H}(\lambda)$ and the Hardy space}\label{sec-4'}

We begin this section with the following concepts.

\bdefe\label{def4'.1}
For $0<p<\infty$, the \textit{Hardy space} $H^p$ is the set of all functions $f$ analytic in $\mathbb{D}$
for which
$$M_p(r,f)=\left\{ \frac{1}{2\pi}\int_{0}^{2\pi}|f(re^{i\theta})|^p
\,d\theta \right \}^{1/p}
$$
is bounded on $0<r<1$.

The space $h^p$ consists of all harmonic mappings $f$ in $\mathbb{D}$ for which $M_p(r,f)$ $(0<r<1)$ are bounded (cf. \cite{Du}).
\edefe





For a harmonic mapping $f=h+\overline{g}$ in $\mathbb{D}$, the
Bloch seminorm is given by (cf. Colonna \cite{Co})
$$\|f\|_{\mathscr{B}_{H}}=\sup_{z\in \mathbb{D}}
(1-|z|^2)\big(|h'(z)|+|g'(z)|\big),
$$
and $f$ is called a Bloch mapping when $\|f\|_{\mathscr{B}_{H}}<\infty$.
In the recent years, the class of harmonic Bloch mappings has been studied extensively together with its
higher dimensional analog (see for example, \cite{CPVW2013,CPW2011,CPW2012,Co} and the references therein).


By Theorem \ref{thm3.1}, we have, for $f\in {\mathcal B}_{H}(\lambda)$,
$$ |f(z)|\leq (1+|b_1|)\int_{0}^{|z|}\left (\frac{1+t}{1-t}\right )^{\lambda}\,dt,
$$
which shows that

\bei
\item $f$ is bounded when $\lambda<1$,

\item $f(z)=O(-\log(1-|z|))$ $(|z|\rightarrow 1)$ when $\lambda=1$,
and

\item $f(z)=O((1-|z|)^{1-\lambda})$ $(|z|\rightarrow 1)$ when
$\lambda>1$.
\eei

Let ${\rm BMOA}$ (resp. ${\rm BMOH}$) denote the class of analytic functions
(resp. harmonic mappings) that have bounded mean oscillation on the
unit disk $\mathbb{D}$ (cf. \cite{Au}). In \cite{Ki}, Kim proved the following result for analytic functions.

\begin{Thm}\label{Thm a}\begin{enumerate}
\item If $\lambda<1$, ${\mathcal B}(\lambda)\cap {\mathcal S}\subset
H^{\infty}$,

\item If $\lambda=1$, ${\mathcal B}(\lambda)\cap {\mathcal S}\subset
{\rm BMOA}$,

\item If $\lambda>1$, ${\mathcal B}(\lambda)\cap {\mathcal S}\subset
H^{p}$ for every $0<p<1/(\lambda-1)$.
\end{enumerate}
\end{Thm}

In order to state a generalization of this result for harmonic mappings, we introduce
$${\mathcal S}_{H_k}=\{f=h+\overline{g}\in {\mathcal S}_{H}:\, f  ~\mbox{ is $k$-quasiconformal}\} $$ for $0\leq k<1$.
We now state the analog of Theorem \ref{Thm a}.

\bthm\label{thm4'.0}
\begin{enumerate}
\item If $\lambda<1$, then ${\mathcal B}_{H}(\lambda)\cap {\mathcal S}_{H}\subset
h^{\infty}$.

\item If $\lambda=1$, then ${\mathcal B}_{H}(\lambda)\cap {\mathcal S}_{H}\subset
{\rm BMOH}$.

\item If $\lambda>1$, then ${\mathcal B}_{H}(\lambda)\cap {\mathcal S}_{H_k} \subset
h^{p}$ for every $0<p<1/(\lambda-1)$.
\end{enumerate}
\ethm

For the proof of Theorem \ref{thm4'.0}, we need some preparation.

\blem\label{lem4'.0} If
$f=h+\overline{g}\in {\mathcal B}_{H}(1)$, then $\|f\|_{\mathscr{B}_{H}}\leq 4(1+|b_1|).$
\elem \bpf
For $f=h+\overline{g}\in {\mathcal B}_{H}(1)$, by Theorem \ref{thm2.1},
we have $h+e^{i\theta}g\in {\mathcal B}(1)$ for each $\theta\in [0,
2\pi]$. It follows from \cite[Theorem 2.1]{Ki} that $
\|h+e^{i\theta}g\|_{\mathscr{B}}\leq 4(1+|b_1|) $, which implies
that $ \|f\|_{\mathscr{B}_{H}}\leq 4(1+|b_1|).$ \epf

In the next lemma, we shall consider the problem of how the integral means
of $k$-quasiconformal harmonic univalent mappings $f$ behaves. Here the integral means of
$f$ is defined by
\be\label{eq-ext1}
I_p(r)=I_p(r,f)=\frac{1}{2\pi}\int_{0}^{2\pi}|f(re^{i\theta})|^p\, d\theta.
\ee
The following lemma is regarded as a generalization of \cite[Proposition 8.1]{Ch1}
to the case of harmonic mappings.

\blem\label{lem4'.1} Let $f\in {\mathcal S}_{H_k}$ and $p>0$. Then
$$I_p(r)\leq \frac{2(1+k^2)(|p-2|+1)}{1-k^2}\int_{0}^{r}M(\rho)^p
\rho^{-1}\, d\rho\;\;\; (0\leq r<1),
$$
where
$$ M(r):=M(r,f)=\max_{0\leq \theta\leq 2\pi}|f(re^{i\theta})|.
$$
\elem \bpf
Let $f=h+\overline{g}\in {\mathcal S}_{H_k}$ and write $z=re^{i\theta}$, where $0\leq r<1$. Writing
$$ |f(z)|^{p}=\left [\big (h(z)+\overline{g(z)}\big  )\big (\overline{h(z)} +g(z)\big  )\right ]^{p/2},
$$
elementary computations give
$$r\frac{\partial}{\partial r}\big(|f(z)|^p\big )=p|f(z)|^{p-2}{\rm
Re}\,\big \{\big(zh'(z)+\overline{zg'(z)}\,\big)\overline{f(z)}\,\big \}
$$
and
$$\frac{\partial}{\partial \theta}\big(|f(z)|^p\big ) =p|f(z)|^{p-2}{\rm Re}\,
\big \{\big(iz h'(z)+\overline{izg'(z)}\,\big)\overline{f(z)}\,\big \}.
$$
Further computations yield
%
\beqq
\Big(r\frac{\partial}{\partial r}\Big)^2\big(|f(z) |^p\big)
&=&p(p-2)|f(z)|^{p-4}\Big({\rm
Re\,}\big \{\big(zh'(z)+\overline{zg'(z)}\,\big) \overline{f(z)}\big\}\Big)^2\\
&&+p|f(z)|^{p-2} {\rm Re\,}\big \{\big(z^2h''(z)+\overline{z^2g''(z)}
+zh'(z)+\overline{zg'(z)}\, \big)\overline{f(z)}\, \big\}\\
&&+p|f(z)|^{p-2}\big|zh'(z)+\overline{zg'(z)}\big|^2
\eeqq
and
%
%
\beqq
\Big(\frac{\partial}{\partial\theta}\Big)^2\big(|f(z)|^p\big)
&=&p(p-2)|f(z)|^{p-4}\Big({\rm Re\,}\big \{\big(iz h'(z)+\overline{i zg'(z)}\,
\big)\overline{f(z)}\,\}\Big)^2\\
&&+p|f(z)|^{p-2}
{\rm Re\,}\big \{\big(-z^2h''(z)-\overline{z^2g''(z)}-zh'(z)-\overline{zg'(z)}\big)
\overline{f(z)}\, \big \}\\
&&+p|f(z)|^{p-2}\big|izh'(z)+\overline{izg'(z)}\big|^2.
\eeqq
Adding the last two expressions shows that

\vspace{8pt}

$\ds \Big(r\frac{\partial}{\partial r}\Big)^2\big(|f(z)|^p\big)
+\Big(\frac{\partial}{\partial\theta}\Big)^2\big(|f(z)|^p\big)
$

\beqq
&=&p(p-2)|f(z)|^{p-4}\left [\Big({\rm
Re\,}\big \{\big( zh'(z)+\overline{zg'(z)}\, \big)\overline{f(z)}\,\big \}\Big)^2\right . \\
&& +\left . \Big({\rm Re\,}\big \{\big(iz h'(z)+\overline{izg'(z)}\, \big)\overline{f(z)}\,\big \}\Big)^2 \right ]\\
&&+p\big|f(z)|^{p-2} \left [\big|zh'(z)+\overline{zg'(z)}\big|^2
+\big|izh'(z)+\overline{izg'(z)}\big|^2\right ]\\
&\leq&2p(1+k^2)(|p-2|+1)r^2|f(z)|^{p-2}|h'(z)|^2.
\eeqq
It follows that

\vspace{8pt}

$\ds \frac{1}{2\pi}\int_{0}^{2\pi}\left [\Big(r\frac{\partial}{\partial
r}\Big)^2\big(|f(re^{i\theta})|^p\big)
+\Big(\frac{\partial}{\partial\theta}\Big)^2\big(|f(re^{i\theta})|^p\big)\right ]d\theta\
$
\beqq
&=&\frac{1}{2\pi}\int_{0}^{2\pi}\Big(r\frac{\partial}{\partial r}\Big)^2\big(|f(re^{i\theta})|^p\big)\,d\theta\\
&=&r\frac{d}{dr}\big(r I_p'(r)\big)\\
&\leq & p(1+k^2)(|p-2|+1)r^2\frac{1}{\pi} \int_{0}^{2\pi}|f(re^{i\theta})|^{p-2}|h'(re^{i\theta})|^2
\,d\theta,
\eeqq
where we have used the fact that the integral corresponding to the second term above vanishes because of the periodicity of the
function involved in the integrand. As $J_f(z)=|h'(z)|^2-|g'(z)|^2$ and
$$|h'(z)|^2 =\frac{|h'(z)|^2}{|h'(z)|^2-|g'(z)|^2}J_f(z)\leq\frac{1}{1-k^2}J_f(z),
$$
we may integrate the last expression on both sides and obtain the inequality
\beqq
r I_p'(r)&\leq & p(1+k^2)(|p-2|+1)\iint_{|z|\leq r}|f(z)|^{p-2}|h'(z)|^2 \,d\sigma (z)\\
&\leq & \frac{p(1+k^2)(|p-2|+1)}{1-k^2}\iint_{|z|\leq  r}|f(z)|^{p-2}J_f(z)\, d\sigma (z),
\eeqq
where $d\sigma (z)=(1/\pi)dx\,dy$ denotes the normalized area element.
Now, we substitute $w=f(z)$. Since $f$ is univalent in $\ID$ and $M(r)=\max_{0\leq \theta\leq
2\pi}|f(re^{i\theta})|$, the last inequality gives
\beqq
r I_p'(r)&\leq & \frac{p(1+k^2)(|p-2|+1)}{1-k^2}\iint_{|w|\leq M(r)} |w|^{p-2}\, d\sigma (w)\\
&=&\frac{2p(1+k^2)(|p-2|+1)}{1-k^2}\int_{0}^{M(r)}t^{p-1}\,dt\\
&=&\frac{2(1+k^2)(|p-2|+1)}{1-k^2}M(r)^p,
\eeqq
which upon integration on both sides shows that
$$ I_p(r)\leq \frac{2(1+k^2)(|p-2|+1)}{1-k^2}\int_{0}^{r}M(\rho)^p
\rho^{-1}\,d\rho.
$$
The desired conclusion follows. \epf


\subsection{The proof of Theorem \ref{thm4'.0}}
Let $f\in {\mathcal B}_{H}(\lambda)$ for some $\lambda<1$. Then, by Corollary \ref{cor3.1}, $f$ is bounded.

Next we assume that $f=h+\overline{g}\in {\mathcal B}_{H}(1)\cap {\mathcal S}_{H}$. Then,
by Lemma \ref{lem4'.0}, it follows that $f$ is Bloch and thus,
$h$ is Bloch, since, for $f=h+\overline{g}\in  {\mathcal S}_{H}$,  $h$ is Bloch if and only if $h$ is BMOA if and
only if $f$ is BMOH (cf. \cite{Au}). Consequently, $f\in {\rm BMOH}$.

Finally, we assume that $f\in {\mathcal B}_{H}(\lambda)\cap {\mathcal S}_{H_k}$ for some $\lambda>1$. Then, by Theorem
\ref{thm3.1}, we deduce that $f(z)=O((1-|z|)^{1-\lambda})$ and thus,
$M(r)=O((1-|z|)^{1-\lambda})$. Furthermore, using Lemma \ref{lem4'.1}, we find that
\beqq
I_p(r)&\leq& \frac{2(1+k^2)(|p-2|+1)}{1-k^2}\int_{0}^{r}M(\rho)^p \rho^{-1}\, d\rho\\
&\leq& \frac{2C (1+k^2)(|p-2|+1)}{1-k^2}\int_{0}^{r}\frac{1}{(1-|z|)^{p(\lambda-1)}} \rho^{-1} \,d\rho
\eeqq
for some positive constant $C$. Hence, $f\in h^p$ if $0<p<1/(\lambda-1)$.
\hfill $\Box$

\bigskip

Obviously, the assertion (1) in Theorem \Ref{Thm a} remains valid if
we replace ${\mathcal B}(\lambda)$ by ${\mathcal B}_{H}(\lambda)$.
By \cite[Theorem 1]{KS1}, we see that the assertion (3) does not hold for ${\mathcal B}_{H}(\lambda)$. On the other hand,
we will show that ${\mathcal B}_{H}(\lambda)$ is contained in some Hardy space.

\bthm\label{thm4'.1}
Let $\lambda\geq 1$. Then ${\mathcal B}_{H}(\lambda)\subset h^p$ with
$0<p<\frac{1}{\lambda^2-1}$.
\ethm

In the above, the expression $\frac{1}{\lambda^2-1}$ is interpreted as $\infty$ when
$\lambda=1$.

\bpf Assume that $f=h+\overline{g}\in {\mathcal B}_{H}(\lambda)$. By
Theorem \ref{thm2.1} and \cite[Theorem 2]{KS1}, for each $\theta$,
$h+e^{i\theta}g\in H^p$ with $0<p<\frac{1}{\lambda^2-1}$. It follows
that $h-g\in H^p$ and $h+g\in H^p$ which implies that $f\in h^p$.
\epf

\bcor
A uniformly locally univalent harmonic mapping $f$ in
$\mathbb{D}$ is contained in the Hardy space $h^p$ for some
$p=p(f)>0$.
\ecor

In \cite{Ki}, Kim also conjectured that the assertion (3) in Theorem
\Ref{Thm a} holds for ${\mathcal B}(\lambda)$.

\section{Coefficient estimates for the class ${\mathcal B}_{H}(\lambda)$}\label{sec-4}

Let
$f(z)=h(z)+\overline{g(z)}=\sum_{n=1}^{\infty}a_nz^n+\sum_{n=1}^{\infty}\overline{b}_n\overline{z}^n$
with $a_1=1$ and $b_1=0$. If $f\in {\mathcal B}_{H}(\lambda)$, then by Theorem \ref{thm2.1}, for each $\theta\in [0,
2\pi]$,
$$\left |\frac{h''(0)+e^{i\theta}g''(0)}{h'(0)+e^{i\theta}g'(0)}\right |\leq 2\lambda,
$$
which shows that $\big | |a_2|-|b_2| \big |\leq |a_2|+|b_2|\leq \lambda$. Of course,
this estimate is sharp because equality holds for $H_\lambda$.

In order to estimate the coefficients of harmonic mappings $f$ in
$\mathbb{D}$, we consider the integral mean $I_p(r,f)$ of $f$ defined by \eqref{eq-ext1},
where $p$ is a real number. For $f(z)=\sum_{n=1}^{\infty}a_nz^n+\sum_{n=1}^{\infty}\overline{b}_n\overline{z}^n
 \in {\mathcal B}_{H}(\lambda)$ with $a_1=1$ and $b_1=0$, by Theorem \ref{thm2.1},
\cite[Theorem 2.3]{KS} and similar arguments as in
\cite[p.~190]{KS}, we have $|a_n+e^{i\theta}b_n|=O\big
(n^{\lambda-1} \big )$ uniformly for $\theta\in [0, 2\pi]$ as
$n\rightarrow \infty$, and then $|a_n|+|b_n|=O(n^{\lambda-1})$ as
$n\rightarrow \infty$. Moreover, if $\lambda<1$ and
$f=h+\overline{g}$ is univalent, then, by Theorem \ref{thm2.1} and
\cite[Corollary 2.4]{KS}, $H_\theta=h+e^{i\theta}g$ is uniformly
bounded for  $\theta\in [0,2\pi]$. So
$${\rm Area}\,(H_\theta(\mathbb{D}))=\pi \Big (1+\sum_{n=2}^{\infty}n|a_n+e^{i\theta}b_n|^2\Big )<\infty,
$$
which implies that $|a_n|+|b_n|=o\big (n^{-1/2}\big )$ as
$n\rightarrow \infty$.

In the following theorem, we improve the exponents in these order
estimates.

\bthm\label{thm4.1} Let
$f(z)=h(z)+\overline{g(z)}=\sum_{n=1}^{\infty}a_nz^n
+\sum_{n=1}^{\infty}\overline{b}_n\overline{z}^n\in {\mathcal
B}_{H}(\lambda)$ with $a_1=1$ and $b_1=0$. Then, for each
$\varepsilon>0$, a real number $p$ and uniformly for $\theta\in [0,
2\pi]$, we have
$$I_p(r, h'+e^{i\theta}g')=O\big ((1-r)^{\alpha(|p|\lambda)-\varepsilon}\big ),
$$
and thus,
$$I_p(r, f)=O\big ((1-r)^{-\alpha(|p|\lambda)-\varepsilon}\big ),~~
|a_n|+|b_n|=O\big (n^{\alpha(\lambda)-1+\varepsilon}\big ),
$$
where $\alpha(\lambda)=\frac{\sqrt{1+4\lambda^2}-1}{2}$.
\ethm

We can prove this theorem by using Theorem \ref{thm2.1} and
\cite[Theroem 3.1]{KS}. Here we omit its detail.

Given a harmonic mapping
$f(z)=h(z)+\overline{g(z)}=\sum_{n=1}^{\infty}a_nz^n+\sum_{n=1}^{\infty}\overline{b}_n\overline{z}^n$
with $a_1=1$ and $b_1=0$ in $\mathbb{D}$, let $\gamma(f)$ denote the
infimum of exponents $\gamma$ such that $|a_n|+|b_n|=O\big
(n^{\gamma-1}\big )$ as $n\rightarrow \infty$, that is,
$$\gamma(f)=\overline{\lim_{n\rightarrow\infty}}\frac{\log n (|a_n|+|b_n|)}{\log n}.
$$

For the subset $X$ of ${\mathcal H}$, we let $\gamma(X)=\sup_{f\in
X}\gamma(f)$. As for the class ${\mathcal S}_b$ of all normalized
bounded univalent functions in $\mathbb{D}$, it is proved that
$0.24<\gamma({\mathcal S}_b)<0.4886$ (cf. \cite{Ca, Ma}), and
conjectured by Carleson and Jones \cite{Ca} that $\gamma({\mathcal
S}_b)=0.25$.  For a  bounded and univalent function $f$,  we note
that the growth of coefficients seems to involve the irregularity of
boundary of image under $f$ (cf. \cite[Chapter 10]{Ch1}), and
Makarov and Pommerenke observed a remarkable phenomenon of phase
transition of the functional $\gamma(f)$ with respect to the
Minkowski dimension of the boundary curve (cf. \cite{Ma}). Recently,
in \cite{KS}, authors established the boundedness of
$\gamma({\mathcal B}(\lambda))$ in terms of $\lambda$. As a
generalization, we consider the class ${\mathcal B}_{H}(\lambda)$
and prove that  $\gamma({\mathcal B}_{H}(\lambda))$ have the same
bound with $\gamma({\mathcal B}(\lambda))$.

 For the class ${\mathcal B}_{H}(\lambda)$, Theorem \ref{thm4.1} implies that
$\gamma({\mathcal B}_{H}(\lambda))\leq \alpha(\lambda)$.
The extremal function $H_\lambda$ satisfies the relation
$\gamma(H_\lambda)=\lambda-1$. By \cite[Example 3.1]{KS}, it follows that
$\gamma({\mathcal B}_{H}(\lambda))\geq 0$ for $\lambda>0$.
Hence, we have

\bthm\label{thm4.2}
For each $\lambda\in (0,\infty)$, we have
$$\max\{\lambda-1,0\}\leq\gamma({\mathcal B}_{H}(\lambda))\leq \alpha(\lambda),
$$
where
$\alpha(\lambda)=\frac{\sqrt{1+4\lambda^2}-1}{2}$. In particular,
$\gamma({\mathcal B}_{H}(\lambda))=O(\lambda^2)$ as
$\lambda\rightarrow 0$.
\ethm

Now we mention a connection with integral means for univalent
analytic functions. For a univalent harmonic mapping $f\in {\mathcal S}_{H}$ and a real number $p$, we let
$$\beta_{f,\theta}(p)=\overline{\lim_{r\rightarrow 1^{-}}}\frac{\log
I_p(r, h'+e^{i\theta}g')}{\log \frac{1}{1-r}}.
$$
Clearly, for an univalent analytic function $f\in {\mathcal S}$,
$$ \beta_{f}(p)=\overline{\lim_{r\rightarrow 1^{-}}}
\frac{\log I_p(r, f')}{\log \frac{1}{1-r}}.
$$
Brennan conjectured that $\beta_f(-2)\leq 1$ for univalent analytic functions $f$ (cf. \cite[Charpter 8]{Ch1}).

As a corollary of Theorem \ref{thm4.1}, we have

\bthm\label{thm4.3}
For $f\in {\mathcal B}_{H}(\lambda)$ and a real number $p$,
$$\beta_{f, \, \theta}(p)\leq \alpha(|p|\lambda)=\frac{\sqrt{1+4p^2\lambda^2}-1}{2}
$$
holds for each $\theta\in [0, 2\pi]$. In particular, the Brennan
conjecture is true for  univalent functions $f$ with $\|T_f\|\leq
\sqrt{2}$.
\ethm

\subsection*{Acknowledgements}
The work of Ms. Jinjing Qiao was supported by National Natural
Science Foundation of China (No. 11501159) and was partially
supported by ``INSA JRD-TATA Fellowship" of the Centre for
International Co-operation in Science (CICS).
The second author is on leave from  IIT Madras.

\end{document}